\documentstyle[amscd, leqno]{amsart}
\theoremstyle{plain}
\newtheorem{theorem}{Theorem}[section]
\newtheorem{proposition}[theorem]{Proposition}
\newtheorem{lemma}[theorem]{Lemma}

\theoremstyle{definition}
\newtheorem{definition}[theorem]{Definition}

\theoremstyle{remark}
\newtheorem{remark}[theorem]{Remark}
\makeatletter

\@addtoreset{equation}{section}
\makeatother

\begin{document}
\title[]
{$K3$ surfaces with involution, 
equivariant analytic torsion, 
and automorphic forms on the moduli space III:
the case $r(M)\geq18$}
\author{Ken-Ichi Yoshikawa}
\address{
Department of Mathematics,
Faculty of Science,
Kyoto University,
Kyoto 606-8502, JAPAN}
\email{yosikawa@@math.kyoto-u.ac.jp}
\address{Korea Institute for Advanced Study,
Hoegiro 87, Dongdaemun-gu,
Seoul 130-722, KOREA}

\thanks{The author is partially supported by the Grants-in-Aid 
for Scientific Research (B) 19340016, JSPS}

\begin{abstract}
We prove the automorphic property of the invariant of $K3$ surfaces with involution,
which we obtained using equivariant analytic torsion, in the case where
the dimension of the moduli space is less than or equal to $2$.
\end{abstract}

\maketitle


\section
{Introduction}
\label{sect:1}
\par
Let $(X,\iota)$ be a $K3$ surface with anti-symplectic 
holomorphic involution and let $H^{2}_{+}(X,{\bf Z})$ be 
the invariant sublattice of $H^{2}(X,{\bf Z})$ with respect to 
the $\iota$-action. By Nikulin \cite{Nikulin83}, 
the topological type of $\iota$ is determined by the isometry 
class of $H^{2}_{+}(X,{\bf Z})$. 
Let $M$ be a sublattice of the $K3$-lattice
and let $M^{\perp}$ be the orthogonal complement of $M$ 
in the $K3$-lattice. 
The pair $(X,\iota)$ is called a $2$-elementary $K3$ surface
of type $M$ if $H^{2}_{+}(X,{\bf Z})$ is isometric to $M$.
In this case, $M$ is a primitive, $2$-elementary, Lorentzian
sublattice of the $K3$-lattice by \cite{Nikulin80}.
Let ${\mathcal M}_{M^{\perp}}^{o}$ be the coarse moduli space 
of $2$-elementary $K3$ surfaces of type $M$. 
By the global Torelli theorem for $K3$ surfaces, the period map
gives an identification between ${\mathcal M}_{M^{\perp}}^{o}$ 
and a Zariski open subset of the modular variety 
$\Omega_{M^{\perp}}^{+}/O^{+}(M^{\perp})$. Here 
$\Omega_{M^{\perp}}^{+}$ is the period domain for 
$2$-elementary $K3$ surfaces of type $M$, which is isomorphic 
to a symmetric bounded domain of type IV of dimension $20-r(M)$,
and $O^{+}(M^{\perp})\subset O(M^{\perp}\otimes{\bf R})$ 
is a certain arithmetic subgroup.
\par
In \cite{Yoshikawa04}, we introduced a real-valued invariant 
$\tau_{M}(X,\iota)$ of $(X,\iota)$, which we obtained using
equivariant analytic torsion \cite{Bismut95} and 
a Bott--Chern secondary class \cite{BGS88}. (See Sect.\ref{sect:2}.) 
Then $\tau_{M}$ gives rise to a function on the coarse moduli 
space ${\mathcal M}_{M^{\perp}}^{o}$. 
\par
Let $r(M)$ be the rank of $M$.
When $r(M)\leq17$, the function $\tau_{M}$
on ${\mathcal M}_{M}^{o}$ is expressed as the Petersson norm
of an automorphic form on $\Omega_{M^{\perp}}^{+}$ 
characterizing the discriminant locus \cite{Yoshikawa04},
where the automorphic form takes its values in a certain 
$O^{+}(M^{\perp})$-equivariant line bundle on 
$\Omega_{M^{\perp}}^{+}$.
The purpose of this note is to extend the automorphic 
property of $\tau_{M}$ to the case $r(M)\geq18$.
\par 
Let $X^{\iota}$ be the set of fixed points of $\iota\colon X\to X$.
If $r(M)\geq18$, $X^{\iota}$ is the disjoint union of finitely
many compact Riemann surfaces, whose total genus is determined
by $M$ (cf. \cite{Nikulin83}). 
Let $g(M)$ be the total genus of $X^{\iota}$.
Then our main result is stated as follows.

\begin{theorem}[Theorem 5.3]\label{Theorem1.1}
There exist an integer $\nu\in{\bf Z}_{>0}$,
an (possibly meromorphic) automorphic form $\Psi_{M}$ 
on $\Omega_{M^{\perp}}^{+}$ of weight $\nu(r(M)-6)$ and 
a Siegel modular form $S_{M}$ on the Siegel upper half space
${\frak S}_{g(M)}$ of weight $4\nu$ such that, 
for every $2$-elementary $K3$ surface $(X,\iota)$ of type $M$, 
$$
\tau_{M}(X,\iota)
=
\|\Psi_{M}(\overline{\varpi}_{M}(X,\iota))\|^{-1/2\nu}
\|S_{M}(\varOmega(X^{\iota})\|^{-1/2\nu}.
$$
Here 
$\overline{\varpi}_{M}(X,\iota)\in{\mathcal M}_{M^{\perp}}^{o}$
denotes the period of $(X,\iota)$, 
$\varOmega(X^{\iota})\in{\frak S}_{g(M)}/Sp_{2g(M)}({\bf Z})$ 
denotes the period of $X^{\iota}$, and $\|\cdot\|$ denotes
the Petersson norm.
\end{theorem}

In \cite{Yoshikawa09a}, we shall use Theorem~\ref{Theorem1.1} 
to give explicit formulae for $\Psi_{M}$ and $S_{M}$. In fact,
$\Psi_{M}$ is expressed as an explicit Borcherds lift
of a certain elliptic modular form and $S_{M}$ is 
expressed as the product of all even theta constants.
\par
This note is organized as follows.
In Sect.\ref{sect:2}, we recall the invariant $\tau_{M}$.
In Sect.\ref{sect:3}, we recall the moduli space of $2$-elementary
$K3$ surfaces of type $M$ and prove a technical result.
In Sect.\ref{sect:4}, we study the singularity of $\tau_{M}$.
In Sect.\ref{sect:5}, we prove Theorem~\ref{Theorem1.1}.
In Sect.\ref{sect:6}, we prove a technical result used 
in the proof of the main theorem for a certain $M$.

\section
{$K3$ surfaces with involution and the invariant $\tau_{M}$}
\label{sect:2}
\par
Let $X$ be a $K3$ surface and let $\iota\colon X\to X$
be a holomorphic involution acting non-trivially on
holomorphic $2$-forms on $X$. The pair $(X,\iota)$ is called 
a $2$-elementary $K3$ surface.
Let ${\Bbb L}_{K3}$ be a fixed even unimodular lattice of 
signature $(3,19)$, which is called a $K3$-lattice. 
Then $H^{2}(X,{\bf Z})$ equipped with the cup-product pairing 
is isometric to ${\Bbb L}_{K3}$.
Let $M\subset{\Bbb L}_{K3}$ be a sublattice.
The pair $(X,\iota)$ is of type $M$ if
the invariant part of $H^{2}(X,{\bf Z})$ with respect
to the $\iota$-action is isometric to $M$.
By \cite{Nikulin80}, there exists a $2$-elementary $K3$ surface
of type $M$ if and only if $M$ is a primitive, $2$-elementary, 
Lorentzian sublattice of ${\Bbb L}_{K3}$. 
\par
Let $(X,\iota)$ be a $2$-elementary $K3$ surface 
of type $M$. Identify ${\bf Z}_{2}$ with 
the subgroup of ${\rm Aut}(X)$ generated by $\iota$. 
Let $\kappa$ be a ${\bf Z}_{2}$-invariant K\"ahler form on $X$. 
Let $\tau_{{\bf Z}_{2}}(X,\kappa)(\iota)$ be the equivariant 
analytic torsion of the trivial Hermitian line bundle on 
$(X,\kappa)$. For the definition and the basic properties of
(equivariant) analytic torsion, we refer the reader to
\cite{RaySinger73}, \cite{BGS88}, \cite{Bismut95},
\cite{KohlerRoessler01}, \cite{Ma00}.
Set ${\rm vol}(X,\kappa):=(2\pi)^{-2}\int_{X}\kappa^{2}/2!$.
Let $\eta$ be a nowhere vanishing holomorphic $2$-form on $X$. 
The $L^{2}$-norm of $\eta$ is defined as 
$\|\eta\|_{L^{2}}^{2}:=(2\pi)^{-2}\int_{X}\eta\wedge\bar{\eta}$.
\par
Let $X^{\iota}:=\{x\in X;\,\iota(x)=x\}$ be the set of fixed 
points of $\iota$ and
let $X^{\iota}=\sum_{i}C_{i}$ be the decomposition
into the connected components. 
By \cite{Nikulin83}, the total genus $g(X^{\iota})$ of $X^{\iota}$ 
depends only on $M$ and hence is denoted by $g(M)$. Set 
${\rm vol}(C_{i},\kappa|_{C_{i}}):=
(2\pi)^{-1}\int_{C_{i}}\kappa|_{C_{i}}$. Let 
$c_{1}(C_{i},\kappa|_{C_{i}})$ be the Chern form of
$(TC_{i},\kappa|_{C_{i}})$ and let
$\tau(C_{i},\kappa|_{C_{i}})$ be the analytic torsion
of the trivial Hermitian line bundle on 
$(C_{i},\kappa|_{C_{i}})$. 
\par
By \cite[Th.\,5.7]{Yoshikawa04}, the real number
$$
\begin{aligned}
\tau_{M}(X,\iota)
&:=
{\rm vol}(X,\kappa)^{\frac{14-r(M)}{4}}
\tau_{{\bf Z}_{2}}(X,\kappa)(\iota)
\prod_{i}{\rm Vol}(C_{i},\kappa|_{C_{i}})
\tau(C_{i},\kappa|_{C_{i}})
\\
&\quad
\times
\exp\left[
\frac{1}{8}\int_{X^{\iota}}
\log\left.\left(\frac{\eta\wedge\bar{\eta}}
{\kappa^{2}/2!}\cdot
\frac{{\rm Vol}(X,\kappa)}{\|\eta\|_{L^{2}}^{2}}
\right)\right|_{X^{\iota}}
c_{1}(X^{\iota},\kappa|_{X^{\iota}})
\right],
\end{aligned}
$$
is independent of the choice of $\kappa$. 
Hence $\tau_{M}(X,\iota)$
is a real-valued invariant of $(X,\iota)$.
We regard $\tau_{M}$ as a function on the moduli space
of $2$-elementary $K3$ surfaces of type $M$.

\section
{The moduli space of $2$-elementary $K3$ surfaces}
\label{sect:3}
\par

\subsection
{The moduli space of $2$-elementary $K3$ surfaces}
\label{subsect:3.1}
\par
By the global Torelli theorem for $K3$ surfaces, 
the period domain for $2$-elementary $K3$ surfaces of 
type $M$ is given by the set
$$
\Omega_{M^{\perp}}
:=
\{[\eta]\in{\bf P}(M^{\perp}\otimes{\bf C});\,
\langle\eta,\eta\rangle=0,\quad
\langle\eta,\overline{\eta}\rangle>0\},
$$
which consists of two connected components
$\Omega_{M^{\perp}}^{+}$ and $\overline{\Omega_{M^{\perp}}^{+}}$.
Since ${\rm sign}(M^{\perp})=(2,20-r(M))$, 
$\Omega_{M^{\perp}}^{+}$ is isomorphic to a symmetric bounded
domain of type IV of dimension $20-r(M)$. Let $O(M^{\perp})$
be the group of isometries of $M^{\perp}$, which acts projectively
on $\Omega_{M^{\perp}}$. Let $O^{+}(M^{\perp})$ be the subgroup
of $O(M^{\perp})$ of index $2$, which preserves the connected
components of $\Omega_{M^{\perp}}$. We define
$$
{\mathcal M}_{M^{\perp}}
:=
\Omega_{M^{\perp}}^{+}/O^{+}(M^{\perp}).
$$
The Baily--Borel--Satake compactification of
${\mathcal M}_{M^{\perp}}$ is denoted by 
${\mathcal M}_{M^{\perp}}^{*}$, which is a normal projective 
variety of dimension $20-r(M)$ with regular part
$({\mathcal M}_{M^{\perp}}^{*})_{\rm reg}$.
\par
Recall that the discriminant locus of $\Omega_{M^{\perp}}^{+}$
is the divisor defined as
$$
{\mathcal D}_{M^{\perp}}
:=
\bigcup_{d\in\Delta_{M^{\perp}}/\pm1}H_{d},
\qquad
H_{d}:=\{[\eta]\in\Omega_{M^{\perp}}^{+};\,
\langle d,\eta\rangle=0\},
$$
where $\Delta_{M^{\perp}}:=\{d\in M^{\perp};\,
\langle d,d\rangle=-2\}$ is the set of roots of $M^{\perp}$.
Let $\overline{\mathcal D}_{M^{\perp}}$ be the divisor of
${\mathcal M}_{M^{\perp}}^{*}$ defined as the closure of 
the image of ${\mathcal D}_{M^{\perp}}$ by the projection
$\varPi_{M}\colon\Omega_{M^{\perp}}^{+}\to
{\mathcal M}_{M^{\perp}}$. 
By \cite[Th.\,1.8]{Yoshikawa04}, the period map induces 
an isomorphism between the coarse moduli space of $2$-elementary 
$K3$ surfaces of type $M$ and the quasi-projective variety
of dimension $20-r(M)$
$$
{\mathcal M}_{M^{\perp}}^{o}
:=
(\Omega_{M^{\perp}}^{+}\setminus{\mathcal D}_{M^{\perp}})
/O^{+}(M^{\perp})
=
{\mathcal M}_{M^{\perp}}\setminus
\overline{\mathcal D}_{M^{\perp}}.
$$
\par
The boundary locus of ${\mathcal M}_{M^{\perp}}^{*}$
is defined as the subvariety:
$$
{\mathcal B}_{M}:=
{\mathcal M}_{M^{\perp}}^{*}\setminus{\mathcal M}_{M^{\perp}}.
$$
Since $\dim{\mathcal B}_{M}=1$ if $r(M)\geq18$ and
$\dim{\mathcal B}_{M}=0$ if $r(M)=19$, ${\mathcal B}_{M}$ is 
a subvariety of ${\mathcal M}_{M^{\perp}}^{*}$
with codimension greater than or equal to $2$ 
when $r(M)\leq 17$ and is a divisor when $r(M)\geq18$.

\subsection
{One parameter families of $2$-elementary $K3$ surfaces}
\label{subsect:3.2}
\par

We need a modification of \cite[Th.\,2.8]{Yoshikawa04},
which shall be used in Sects.\ref{sect:4} and \ref{sect:6}.

\begin{theorem}\label{Theorem3.1}
Let $C\subset{\mathcal M}_{M^{\perp}}^{*}$ be an irreducible 
projective curve. 
\begin{itemize}
\item[(1)]
There exist a smooth projective curve $B$, 
a morphism $\varphi\colon B\to{\mathcal M}_{M^{\perp}}^{*}$,
an irreducible projective threefold $X$ with an involution 
$\theta\colon X\to X$, and a surjective morphism $f\colon X\to B$
with the following properties:
\begin{itemize}
\item[(a)]
$\varphi(B)=C$.
\item[(b)]
The involution $\theta\colon X\to X$ preserves the fibers of 
$f\colon X\to B$.
\item[(c)]
There is a non-empty Zariski open subset $B^{o}\subset B$ 
such that $(X_{b},\theta|_{X_{b}})$ is a $2$-elementary 
$K3$ surface of type $M$ with period $\varphi(b)$
for $b\in B^{o}$.
\end{itemize}
\item[(2)]
Let $p\colon{\mathcal Z}\to\varDelta$ be a proper surjective
projective morphism from a smooth threefold to the unit disc and
let $\iota\colon{\mathcal Z}\to{\mathcal Z}$ be a holomorphic
involution preserving the fibers $Z_{t}=p^{-1}(t)$ of $p$. 
Assume that $(Z_{t},\iota|_{Z_{t}})$ is a $2$-elementary
$K3$ surface for all $t\in\varDelta^{*}:=\varDelta\setminus\{0\}$ 
and that the period map for 
$p\colon({\mathcal Z},\iota)|_{\varDelta^{*}}\to\varDelta^{*}$
extends to a non-constant holomorphic map 
$\gamma\colon\varDelta\to C$. Let $\nu\in{\bf Z}_{\geq0}$.
Then there exist $\varphi\colon B\to C$, $f\colon X\to B$,
$\theta\colon X\to X$ as above in (1) 
and a point ${\frak p}\in\varphi^{-1}(\gamma(0))$ and
an isomorphism of germs 
$\psi\colon(\varDelta,0)\cong(B,{\frak p})$ 
with the following properties:
\begin{itemize}
\item[(d)]
$(X_{\frak p},\theta|_{X_{\frak p}})\cong(Z_{0},\iota|_{Z_{0}})$
\item[(e)]
The maps of germs 
$\varphi\colon(B,{\frak p})\to(C,\gamma(0))$ 
and $\gamma\colon(\varDelta,0)\to(C,\gamma(0))$ have
the same $\nu$-jets:
For any $F\in{\mathcal O}_{C,\gamma(0)}$,
$$
F\circ\varphi\circ\psi(t)-F\circ\gamma(t)\in 
t^{\nu+1}{\bf C}\{t\}.
$$
\item [(f)]
Let ${\rm Def}(Z_{0})$ be the Kuranishi space of $Z_{0}$ and
let $\mu_{f}\colon(B,{\frak p})\to{\rm Def}(Z_{0})$ and
$\mu_{p}\colon(\varDelta,0)\to{\rm Def}(Z_{0})$
be the maps of germs induced by the deformations
$f\colon(X,X_{\frak p})\to(B,{\frak p})$ and 
$p\colon({\mathcal Z},Z_{0})\to(\varDelta,0)$, respectively.
Then $\mu_{f}$ and $\mu_{p}$ have the same $\nu$-jets:
For any $F\in{\mathcal O}_{{\rm Def}(Z_{0})}$,
$$
F\circ\mu_{f}\circ\psi(t)-F\circ\mu_{p}(t)\in 
t^{\nu+1}{\bf C}\{t\}.
$$
\end{itemize}
\end{itemize}
\end{theorem}

\begin{pf}
We follow \cite[Th.\,2.8]{Yoshikawa04}.
By the same argument as in
\cite[Proof of Th.\,2.8 (Step 1) and Claim 1]{Yoshikawa04},
there exist an irreducible projective variety $T$ and
a family of projective surfaces with involution
$\pi\colon({\mathcal X},{\mathcal I})\to T$ 
with the following properties:
\begin{itemize}
\item[(i)]
Let $D\subset T$ be the discriminant locus of 
$\pi\colon{\mathcal X}\to T$ and define
$T^{o}:=T\setminus({\rm Sing}\,T\cup D)$.
Then $({\mathcal X}_{t},{\mathcal I}_{t})$
is a $2$-elementary $K3$ surface of type $M$ for all $t\in T^{o}$.
\item[(ii)]
Let 
$\overline{\varpi}_{T^{o}}\colon T^{o}\to
{\mathcal M}_{M^{\perp}}^{o}$ be the period map for 
$\pi|_{T^{o}}\colon({\mathcal X}|_{T^{o}},{\mathcal I}|_{T^{o}})
\to T^{o}$. Then $\overline{\varpi}_{T^{o}}(T^{o})\subset C$
and $\overline{\varpi}_{T^{o}}(T^{o})$ contains a non-empty
open subset of $C$. 
\item[(iii)]
The period map $\overline{\varpi}_{T^{o}}\colon T^{o}\to C$ 
extends to a rational map 
$\overline{\varpi}_{T}\colon T\dashrightarrow C$.
\item[(iv)]
In (2), there is a map $c\colon\varDelta\to T$ with
$c(\varDelta^{*})\subset T^{o}$ such that
$p\colon({\mathcal Z},\iota)\to\varDelta$ is induced from 
$\pi\colon({\mathcal X},{\mathcal I})\to T$ by $c$.
\end{itemize}
\par{\bf (1) }
Let $\Gamma\subset T\times C$ be the closure of the graph of
$\overline{\varpi}_{T^{o}}$. Let $B$ be a smooth projective curve
and let $h\colon B\to\Gamma$ be a holomorphic map with
${\rm pr}_{2}(h(B))=C$. We set 
$\varphi:={\rm pr}_{2}\circ h\colon B\to C$.
Let $\pi_{B}\colon({\mathcal X}\times_{T}B,
{\mathcal I}\times{\rm id}_{B})\to B$ be the family of algebraic 
surfaces with involution induced from
$\pi\colon({\mathcal X},{\mathcal I})\to T$ by
${\rm pr}_{1}\circ h\colon B\to T$.
Then the period map for
$\pi_{B}\colon({\mathcal X}\times_{T}B,
{\mathcal I}\times{\rm id}_{B})\to B$ is given by
$\overline{\varpi}_{T}\circ{\rm pr}_{1}\circ h$.
Since $\Gamma\subset T\times C$ is the closure of the graph of
$\overline{\varpi}_{T^{o}}$, we get
$\overline{\varpi}_{T}\circ{\rm pr}_{1}\circ h=
{\rm pr}_{2}\circ h=\varphi$. If we set
$B^{o}:=\nu^{-1}(B\cap(T^{o}\times C))$, then
$\pi_{B}\colon({\mathcal X}\times_{T}B,
{\mathcal I}\times{\rm id}_{B})\to B$ satisfies (a), (b), (c).
This proves (1).
\par{\bf (2) }
To prove (2), we must choose $B$ more carefully as in
\cite[Proof of Th.\,2.8 Claim 2]{Yoshikawa04}.
Let $\sigma\colon\varDelta\to\Gamma$ be the map 
defined as $\sigma(t):=(c(t),\gamma(t))$ for $t\in\varDelta$.
Let $\varSigma\colon\widetilde{\Gamma}\to\Gamma$ be a resolution
such that $\widetilde{\Gamma}$ is projective. Since
$c(\varDelta^{*})\subset T\setminus{\rm Sing}\,T$ and hence 
$\sigma(\varDelta^{*})\subset\Gamma\setminus{\rm Sing}\,\Gamma$,
$\sigma$ lifts to a holomorphic map
$\widetilde{\sigma}\colon\varDelta\to\widetilde{\Gamma}$ 
such that $\sigma=\varSigma\circ\widetilde{\sigma}$.
By \cite[Th.\,1.1]{DLS94}, there exist a pointed smooth 
projective curve $(B_{\nu},{\frak p}_{\nu})$, a holomorphic map 
$\widetilde{h}_{\nu}\colon B_{\nu}\to\widetilde{\Gamma}$ 
and an isomorphism of germs
$\psi\colon(\varDelta,0)\cong(B,{\frak p})$
such that for any 
$G\in{\mathcal O}_{\widetilde{\Gamma},\widetilde{\sigma}(0)}$,
\begin{equation}\label{eqn:3.1}
G\circ\widetilde{h}_{\nu}\circ\psi(t)
-
G\circ\widetilde{\sigma}(t)
\in 
t^{\nu+1}{\bf C}\{t\}.
\end{equation}
\par
Set $B:=B_{\nu}$,
$h:=\varSigma\circ\widetilde{h}_{\nu}\colon B\to\Gamma$ and 
we consider the family of $2$-elementary $K3$ surfaces 
$\pi_{B}\colon({\mathcal X}\times_{T}B,
{\mathcal I}\times{\rm id}_{B})\to B$ of type $M$.
By construction, we get (a), (b), (c), (d).
Since 
\begin{equation}\label{eqn:3.2}
\varphi={\rm pr}_{2}\circ h=
({\rm pr}_{2}\circ\varSigma)\circ\widetilde{h}_{\nu},
\qquad
\gamma={\rm pr}_{2}\circ\sigma=
({\rm pr}_{2}\circ\varSigma)\circ\widetilde{\sigma},
\end{equation}
we get by \eqref{eqn:3.1}, \eqref{eqn:3.2}
\begin{equation}\label{eqn:3.3}
F\circ\varphi\circ\psi(t)-F\circ\gamma(t)
=
(F\circ{\rm pr}_{2}\circ\varSigma)\circ\widetilde{h}(t)
-
(F\circ{\rm pr}_{2}\circ\varSigma)\circ\widetilde{\sigma}(t)
\in
t^{\nu+1}{\bf C}\{t\}.
\end{equation}
This proves (e).
Let $\mu_{\pi}\colon(T,c(0))\to{\rm Def}(Z_{0})$ be the map
induced by the deformation
$\pi\colon({\mathcal X},X_{c(0)})\to(T,c(0))$.
Since
\begin{equation}\label{eqn:3.4}
\mu_{f}=\mu_{\pi}\circ{\rm pr}_{1}\circ h=
(\mu_{\pi}\circ{\rm pr}_{1}\circ\varSigma)\circ\widetilde{h},
\qquad
\mu_{p}=\mu_{\pi}\circ{\rm pr}_{1}\circ\sigma=
(\mu_{\pi}\circ{\rm pr}_{1}\circ\varSigma)\circ\widetilde{\sigma},
\end{equation}
we get by \eqref{eqn:3.1}, \eqref{eqn:3.4}
\begin{equation}\label{eqn:3.5}
F\circ\varphi\circ\psi(t)-F\circ\gamma(t)
=
(F\circ{\rm pr}_{1}\circ\varSigma)\circ\widetilde{h}(t)
-
(F\circ{\rm pr}_{1}\circ\varSigma)\circ\widetilde{\sigma}(t)
\in
t^{\nu+1}{\bf C}\{t\}.
\end{equation}
This proves (f). This completes the proof of (2).
\end{pf}

\section
{The singularity of $\tau_{M}$}
\label{sect:4}
\par
We prove the logarithmic divergence of $\tau_{M}$ for any
one-parameter degeneration of $2$-elementary $K3$ surfaces
of type $M$. For this, we recall the following:

\begin{theorem}\label{Theorem4.1}
Let $\pi\colon X\to S$ be a proper surjective holomorphic map 
from a connected projective algebraic manifold $X$ of 
dimension $n+1$ to a compact Riemann surface $S$. 
Let $G$ be a finite group. 
Assume that $G$ acts holomorphically on $X$ and trivially 
on $S$ and that $\pi\colon X\to S$ is $G$-equivariant. Hence
$G$ preserves all the fibers $X_{s}:=\pi^{-1}(s)$, $s\in S$. 
Let $\Delta:=\{s\in S;\,{\rm Sing}(X_{s})\not=\emptyset\}$ 
be the discriminant locus.
\par
Let $h_{X}$ be a $G$-invariant K\"ahler metric on $X$ and set
$h_{s}:=h_{X}|_{X_{s}}$ for $t\in S\setminus\Delta$. 
Let $\tau_{G}(X_{s},h_{s})(g)$ be the equivariant analytic 
torsion of the trivial Hermitian line bundle on $(X_{s},h_{s})$.
Let $t$ be a local coordinate of $S$ centered at $0\in\Delta$.
If $N$ is the order of $g\in G$, then there exists 
$\beta_{g}(\pi,X_{0})\in\sum_{0\leq k<N}{\bf Q}\exp(2\pi ik/N)$
such that
$$
\log\tau_{G}(X_{t},h_{t})(g)
=
\beta_{g}(\pi,X_{0})\,\log|t|^{2}+O(\log(-\log|t|))
\qquad
(t\to0).
$$
\end{theorem}

\begin{pf}
See \cite[Th.\,1.1 and Cor.\,6.10]{Yoshikawa09b}.
\end{pf}

\begin{theorem}\label{Theorem4.2}
Let $(S,0)$ be a pointed smooth projective curve
equipped with a coordinate neighborhood $(U,t)$
centered at $0$,
let $X$ be a smooth projective threefold equipped
with a holomorphic involution $\theta\colon X\to X$, and 
let $\pi\colon X\to S$ be a surjective holomorphic map.
Assume the following:
\begin{itemize}
\item[(1)]
the projection $\pi\colon X\to S$ is ${\bf Z}_{2}$-equivariant 
with respect to the ${\bf Z}_{2}$-action on $X$ induced by 
$\theta$ and with respect to the trivial ${\bf Z}_{2}$-action 
on $S$.
\item[(2)]
$(X_{t},\theta|_{X_{t}})$ is a $2$-elementary $K3$ surfaces 
of type $M$ for all $t\in U\setminus\{0\}$.
\end{itemize}
Then there exists $\alpha\in{\bf Q}$ such that
$$
\log\tau_{M}(X_{t},\theta|_{X_{t}})
=
\alpha\,\log|t|^{2}+O\left(\log(-\log|t|^{2})\right)
\qquad
(t\to 0).
$$
\end{theorem}

\begin{pf}
Set $\theta_{t}:=\theta|_{X_{t}}$.
Let $h_{X}$ be a ${\bf Z}_{2}$-invariant K\"ahler metric
on $X$ with K\"ahler form $\omega_{X}$
and set $\omega_{t}:=\omega_{X}|_{X_{t}}$.
By Theorem~\ref{Theorem4.1}, 
there exists $\beta\in{\bf Q}$ such that
\begin{equation}\label{eqn:4.1}
\log\tau_{{\bf Z}_{2}}(X_{t},\omega_{t})(\theta_{t})
=
\beta\,\log|t|^{2}+O\left(\log(-\log|t|^{2})\right)
\qquad
(t\to0).
\end{equation}
\par
Let $X^{\theta}$ be the set of fixed points of 
$\theta\colon X\to X$ and let $\Delta\subset S$ be 
the discriminant locus of $\pi\colon X\to S$. 
By the ${\bf Z}_{2}$-equivariance of $\pi$, 
we have the decomposition 
$$
X^{\theta}=X^{\theta}_{H}\amalg X^{\theta}_{V},
$$ 
where $\pi(X^{\theta}_{V})\subset\Delta$ and 
$\pi|_{X_{H}^{\theta}}$ is a surjective 
map from any component of $X^{\theta}_{H}$ to $S$.
Set $Y:=X^{\theta}_{H}$ and $f:=\pi|_{X^{\theta}_{H}}$. 
Then $Y$ is a smooth complex surface and
$f\colon Y\to S$ is a proper surjective holomorphic map
such that $Y_{t}=X_{t}^{\theta_{t}}$ is the disjoint union of 
compact Riemann surfaces for $t\in U\setminus\{0\}$.
It follows from Theorem~\ref{Theorem4.1} again that
there exists $\gamma\in{\bf Q}$ with
\begin{equation}\label{eqn:4.2}
\log\tau(X_{t}^{\theta_{t}},\omega_{t}|_{X_{t}^{\theta_{t}}})
=
\log\tau(Y_{t},\omega_{t}|_{Y_{t}})
=
\gamma\,\log|t|^{2}+O\left(\log(-\log|t|^{2})\right)
\quad
(t\to0).
\end{equation}
\par
Let $K_{X/S}:=
\Omega_{X}^{3}\otimes(\pi^{*}\Omega_{S}^{1})^{-1}$
be the relative canonical bundle. Then the direct image sheaf
$\pi_{*}K_{X/S}$ is locally free on $S$ by e.g.
\cite[Th.\,6.10 (iv)]{Takegoshi95}.
By assumption (2), $\pi_{*}K_{X/S}$ has rank one.
By shrinking $U$ if necessary, there exists
$\Xi\in H^{0}(\pi^{-1}(U),\Omega^{3}_{X})$ such that
$\eta_{X/S}:=\Xi\otimes(\pi^{*}dt)^{-1}$ generates 
$\pi_{*}K_{X/S}$ as an ${\mathcal O}_{S}$-module 
over $U$. In particular, we may assume 
$\eta_{X/S}|_{X_{t}}\not=0$ for $t\not=0$. 
Since $K_{X/S}|_{X_{t}}$ is trivial for $t\not=0$ by (2),
this implies that $\eta_{X/S}|_{X_{t}}$ is nowhere vanishing
on $X_{t}$, $t\not=0$. Hence ${\rm div}(\Xi)\subset X_{0}$.
We set
$\eta_{t}:={\rm Res}_{X_{t}}[\Xi/(\pi-t)]\in 
H^{0}(X_{t},\Omega^{2}_{X_{t}})$ for $t\in U$. Then
$\eta_{X/S}|_{X_{t}}=\Xi\otimes(\pi^{*}dt)^{-1}|_{X_{t}}$
is identified with $\eta_{t}$. 
\par
We prove the existence of $\delta\in{\bf Q}$ such that 
as $t\to 0$
\begin{equation}\label{eqn:4.3}
\int_{X_{t}^{\theta_{t}}}
\log\left.\left(\frac{\eta_{t}\wedge\bar{\eta}_{t}}
{\omega_{t}^{2}/2!}\cdot
\frac{{\rm Vol}(X_{t},\omega_{t})}{\|\eta_{t}\|_{L^{2}}^{2}}
\right)\right|_{X_{t}^{\theta_{t}}}
c_{1}(X_{t}^{\theta_{t}},\omega_{t}|_{X_{t}^{\theta_{t}}})
=
\delta\,\log|t|^{2}+O\left(\log(-\log|t|^{2})\right).
\end{equation}
\par
Let $\Sigma_{\pi}\subset X$ be the critical locus of $\pi$
and let $TX/S:=\ker\pi_{*}|_{X\setminus\Sigma_{\pi}}$ 
be the relative tangent bundle of $\pi\colon X\to S$.
Let $h_{X/S}:=h_{X}|_{TX/S}$ be the Hermitian metric on $TX/S$
induced from $h_{X}$ and let $\omega_{X/S}$ be the $(1,1)$-form 
on $TX/S$ associated to $h_{X/S}$. We identify $\omega_{X/S}$
with the family of K\"ahler forms $\{\omega_{t}\}_{t\in S}$.
Let $N_{X_{t}/X}^{*}$ be the conormal bundle
of $X_{t}$ in $X$ for $t\in U\setminus\{0\}$.
Since $d\pi=\pi^{*}dt\in H^{0}(X_{t},N^{*}_{X_{t}/X})$
generates $N^{*}_{X_{t}/X}$ for $t\in U\setminus\{0\}$, 
$N^{*}_{X_{t}/X}$ is trivial in this case.
Since the Hermitian metric on $\Omega^{1}_{X_{t}}$ is 
induced from $h_{X}$ via the $C^{\infty}$ identification
$\Omega^{1}_{X_{t}}\cong(N_{X_{t}/X}^{*})^{\perp}$
and since $(\omega_{X/S}^{2}/2!)|_{X_{t}}$ is the volume 
form on $X_{t}$, we get on $X\setminus\Sigma_{\pi}$
\begin{equation}\label{eqn:4.4}
\frac{\omega_{X}^{3}}{3!}
=
\frac{\omega_{X/S}^{2}}{2!}\wedge
\left(i\,\frac{d\pi}{\|d\pi\|}\wedge
\frac{\overline{d\pi}}{\|d\pi\|}\right).
\end{equation}
Since $\Xi|_{X_{t}}=\eta_{t}\otimes d\pi$, we get the following 
equation on $X\setminus\Sigma_{\pi}$ by \eqref{eqn:4.4}
\begin{equation}\label{eqn:4.5}
\frac{\eta_{X/S}\wedge\overline{\eta_{X/S}}}
{\omega_{X/S}^{2}/2!}
=
\frac{(-1)^{3}i^{3}\,\Xi\wedge\overline{\Xi}}
{(\omega_{X/S}^{2}/2!)\wedge(i\,d\pi\wedge\overline{d\pi})}
=
\frac{(-1)^{3}i^{3}\,\Xi\wedge\overline{\Xi}}{\omega_{X}^{3}/3!}
\cdot
\frac{1}{\|d\pi\|^{2}}
=
\frac{\|\Xi\|^{2}}{\|d\pi\|^{2}}.
\end{equation}
\par
Let $\Sigma_{f}\subset Y$ be the critical locus of 
$f\colon Y\to S$ and 
let $h_{TY/S}$ be the metric on the relative tangent bundle 
$TY/S:=\ker f_{*}|_{Y\setminus\Sigma_{f}}$ induced from $h_{X}$ 
via the inclusion $TY/S\subset TY\subset TX|_{Y}$. Define
\begin{equation}\label{eqn:4.6}
\begin{aligned}
{\mathcal A}(X/S)
&:=
f_{*}\left[\log\left.\left(
\frac{\eta_{X/S}\wedge\overline{\eta_{X/S}}}
{\omega_{X/S}^{2}/2!}\right)
\right|_{Y\setminus f^{-1}(\pi(\Sigma_{\pi}))}
c_{1}(TY/S,h_{TY/S})
\right]
\\
&\quad
+\chi(Y_{\rm gen})\log
\frac{{\rm Vol}(Y_{\rm gen},\omega_{X}|_{Y_{\rm gen}})}
{\|\eta_{X/S}\|_{L^{2}}^{2}},
\end{aligned}
\end{equation}
where $Y_{\rm gen}$ denotes a general fiber of $f\colon Y\to S$
and $\chi(Y_{\rm gen})$ denotes its topological Euler number.
By \cite[Th.\,6.8]{Yoshikawa09b}, 
there exists $\epsilon_{1}\in{\bf Q}$ such that
\begin{equation}\label{eqn:4.7}
\log\|\eta_{X/S}\|_{L^{2}}^{2}
=
\epsilon_{1}\,\log|t|^{2}+O\left(\log(-\log|t|^{2})\right)
\qquad
(t\to0).
\end{equation}
\par
Let $\varpi\colon{\bf P}(TY)\to Y$ be the projection from
the projective tangent bundle of $Y$ to $Y$.
Let $q\colon\widetilde{Y}\to Y$ be the resolution of 
the indeterminacy of the Gauss map 
$\nu\colon Y\setminus\Sigma_{f}\ni y\to
[T_{y}Y_{f(y)}]\in{\bf P}(TY)$ (cf. \cite[Sect.\,2]{Yoshikawa07})
and set $\widetilde{f}:=f\circ q\colon\widetilde{Y}\to S$ and
$\widetilde{\nu}:=\nu\circ q\colon\widetilde{Y}\to{\bf P}(TY)$. 
Then $\widetilde{f}$ and $\widetilde{\nu}$ are holomorphic maps.
Let ${\mathcal L}\to{\bf P}(TY)$ be the universal line bundle
and let $h_{\mathcal L}$ be the metric on ${\mathcal L}$ 
induced from $\varpi^{*}h_{Y}$ via the inclusion
${\mathcal L}\subset\varpi^{*}TY$. Then
\begin{equation}\label{eqn:4.8}
c_{1}(TY/S,h_{TY/S})
=
\widetilde{\nu}^{*}c_{1}({\mathcal L},h_{\mathcal L}).
\end{equation}
Substituting \eqref{eqn:4.5} into \eqref{eqn:4.6}, we get
\begin{equation}\label{eqn:4.9}
\begin{aligned}
{\mathcal A}(X/S)
&=
f_{*}\left[\log
\left(\frac{\|\Xi\|^{2}}{\|d\pi\|^{2}}
\right)\,
c_{1}(TY/S,h_{TY/S})\right]
-\chi_{\rm top}(Y_{\rm gen})\log\|\eta_{X/S}\|_{L^{2}}^{2}
+O(1)
\\
&=
\widetilde{f}_{*}\left[\log(q^{*}\|\Xi\|^{2})\,
\widetilde{\nu}^{*}c_{1}({\mathcal L},h_{\mathcal L})
\right]
-
\widetilde{f}_{*}\left[\log(q^{*}\|d\pi\|^{2})\,
\widetilde{\nu}^{*}c_{1}({\mathcal L},h_{\mathcal L})
\right]
\\
&\quad
-\epsilon_{1}\,\chi_{\rm top}(Y_{\rm gen})\,\log|t|^{2}
+O\left(\log(-\log|t|^{2})\right),
\end{aligned}
\end{equation}
where we used \eqref{eqn:4.7} and \eqref{eqn:4.8} 
to get the second equality.
Since $q^{*}\Xi$ is a holomorphic section of 
the holomorphic line bundle $q^{*}\Omega_{X}^{3}$ with
${\rm div}(q^{*}\Xi)\subset\widetilde{\pi}^{-1}(0)$, 
there exists by \cite[Lemma 4.4]{Yoshikawa07}
a constant $\epsilon_{2}\in{\bf Q}$ such that
\begin{equation}\label{eqn:4.10}
\widetilde{f}_{*}\left[\log(q^{*}\|\Xi\|^{2})\,
\widetilde{\nu}^{*}c_{1}({\mathcal L},h_{\mathcal L})
\right]
=
\epsilon_{2}\log|t|^{2}+O(1)
\qquad
(t\to 0).
\end{equation}
By \cite[Cor.\,4.6]{Yoshikawa07}, there exists
$\epsilon_{3}\in{\bf Q}$ such that
\begin{equation}\label{eqn:4.11}
\widetilde{f}_{*}\left[\log(q^{*}\|d\pi\|^{2})\,
\widetilde{\nu}^{*}c_{1}({\mathcal L},h_{\mathcal L})
\right]
=
\epsilon_{3}\log|t|^{2}+O(1)
\qquad
(t\to0).
\end{equation}
Setting $\delta:=\epsilon_{2}-\epsilon_{3}-
\epsilon_{1}\,\chi_{\rm top}(Y_{\rm gen})\in{\bf Q}$,
we get \eqref{eqn:4.3} by \eqref{eqn:4.9}, \eqref{eqn:4.10}, 
\eqref{eqn:4.11}.
\par
By the definition of $\tau_{M}$,
the result follows from \eqref{eqn:4.1}, \eqref{eqn:4.2}, 
\eqref{eqn:4.3}.
\end{pf}

\begin{theorem}\label{Theorem4.3}
Let $C\subset{\mathcal M}_{M^{\perp}}^{*}$ be an irreducible 
projective curve intersecting 
$\overline{\mathcal D}_{M^{\perp}}\cup{\mathcal B}_{M}$ properly. 
Let ${\frak b}\in
C\cap(\overline{\mathcal D}_{M^{\perp}}\cup{\mathcal B}_{M})$
and let $C_{\frak b}=\bigcup_{i\in I}C^{(i)}_{\frak b}$
be the irreducible decomposition of the set germ 
$C_{\frak b}=(C,{\frak b})$.
Let $\nu^{(i)}\colon(\varDelta,0)\to C^{(i)}_{\frak b}$ 
be the normalization. Then there exists 
$\alpha_{\frak b}^{(i)}\in{\bf Q}$ such that as $t\to0$,
$$
\log\tau_{M}(\nu^{(i)}(t))
=
\alpha^{(i)}_{\frak b}\,\log|t|+O(\log(-\log|t|)).
$$
\end{theorem}

\begin{pf}
Let $f\colon(X,\theta)\to B$ be the family of $2$-elementary
$K3$ surfaces of type $M$ with period map $\varphi\colon B\to C$
as in Theorem~\ref{Theorem3.1} (1). 
By \cite[Th.\,13.4]{BierstoneMilman97}, there exists 
a resolution of the singularities $\mu\colon\widetilde{X}\to X$ 
such that $\theta$ lifts to an involution 
$\widetilde{\theta}\colon\widetilde{X}\to\widetilde{X}$.
We set $\widetilde{f}:=f\circ\mu$. Since $\mu$ is an isomorphism
outside the singular fibers of $f$, the period map for
$\widetilde{f}\colon(\widetilde{X},\widetilde{\theta})\to B$
coincides with $\varphi\colon B\to C$.
Replacing $f\colon(X,\theta)\to B$ by 
$\widetilde{f}\colon(\widetilde{X},\widetilde{\theta})\to B$
if necessary, we may assume that $X$ is smooth.
\par
For $i\in I$,
let ${\frak p}^{(i)}\in\varphi^{-1}({\frak b})$ be such that
$\varphi(B_{{\frak p}^{(i)}})=C^{(i)}_{\frak b}$.
Let $(V^{(i)},s)$ be a coordinate neighborhood of
${\frak p}^{(i)}$ in $B$ with $s({\frak p}^{(i)})=0$. 
Let $\varphi^{(i)}\colon V^{(i)}\to\varDelta$ be the holomorphic
map such that $\varphi^{(i)}=(\nu^{(i)})^{-1}\circ\varphi$ on 
$V^{(i)}\setminus\{{\frak p}^{(i)}\}$. 
There exists $m_{i}\in{\bf Z}_{>0}$ and 
$\epsilon_{i}(s)\in{\bf C}\{s\}$
such that $t\circ\varphi^{(i)}(s)=s^{m_{i}}\epsilon_{i}(s)$ 
and $\epsilon_{i}(0)\not=0$.
By Theorem~\ref{Theorem4.2} applied to the family 
$f\colon(X,\theta)\to B$, 
there exists $\alpha_{i}\in{\bf Q}$ such that
$$
\log\tau_{M}(\nu^{(i)}\circ\varphi^{(i)}(s))
=
\alpha_{i}\,\log|s|+O(\log(-\log|s|))
\qquad
(s\to0).
$$
This, together with the relation 
$t\circ\varphi^{(i)}(s)=s^{m_{i}}\epsilon_{i}(s)$, 
yields the desired estimate with 
$\alpha^{(i)}_{\frak b}=\alpha_{i}/m_{i}$.
\end{pf}

\section
{The automorphic property of $\tau_{M}$: the case $r(M)\geq18$}
\label{sect:5}
\par
In \cite[Main Th.]{Yoshikawa04}, we proved that $\tau_{M}$ is
expressed as the Petersson norm of an automorphic form
on the period domain for $2$-elementary $K3$ surfaces of
type $M$ if $r(M)\leq17$. In this section, we extend
this result when $r(M)\geq18$. For $n\in{\bf Z}$,
$\langle n\rangle$ denotes the $1$-dimensional lattice ${\bf Z}$
equipped with the bilinear form $\langle x,y\rangle=nxy$.
We denote by ${\Bbb U}$ the $2$-dimensional lattice associated
to the matrix $\binom{0\,1}{1\,0}$.

\subsection
{Automorphic forms on the moduli space}
\label{subsect:5.1}
\par
We fix a vector $l_{M^{\perp}}\in M^{\perp}\otimes{\bf R}$ 
with $\langle l_{M^{\perp}},l_{M^{\perp}}\rangle\geq0$ and set
$$
j_{M^{\perp}}(\gamma,[\eta])
:=
\frac{\langle\gamma(\eta),l_{M^{\perp}}\rangle}
{\langle\eta,l_{M^{\perp}}\rangle}
\qquad
[\eta]\in\Omega_{M^{\perp}}^{+},
\quad
\gamma\in O^{+}(M^{\perp}).
$$
Since $\langle l_{M^{\perp}},l_{M^{\perp}}\rangle\geq0$,
$j_{M^{\perp}}(\gamma,\cdot)$ is a nowhere vanishing holomorphic
function on $\Omega_{M^{\perp}}^{+}$.

\begin{definition}\label{Definition5.1}
A holomorphic function 
$F\in{\mathcal O}(\Omega_{M^{\perp}}^{+})$ is an automorphic form 
on $\Omega_{M^{\perp}}^{+}$ for $O^{+}(M^{\perp})$ of weight $\nu$ 
if the following two conditions are satisfied:
\begin{itemize}
\item[(i)]
There exists a unitary character
$\chi\colon O^{+}(M^{\perp})\to U(1)$ such that 
$$
F([\gamma(\eta)])
=
\chi(\gamma)\,j_{M^{\perp}}(\gamma,[\eta])^{\nu}\,F([\eta]),
\qquad
[\eta]\in\Omega_{M^{\perp}}^{+},
\quad
\gamma\in O^{+}(M^{\perp}).
$$
\item[(ii)]
Denote by $\|F\|^{2}\in C^{\infty}(\Omega_{M^{\perp}}^{+})$ 
the Petersson norm of $F$, which is regarded as 
a $C^{\infty}$ function on ${\mathcal M}_{M^{\perp}}$ 
in the sense of orbifolds. Then 
$\log\|F\|^{2}\in 
L^{1}_{\rm loc}(({\mathcal M}_{M^{\perp}}^{*})_{\rm reg})$
and there exists an effective divisor $D$ on
${\mathcal M}_{M^{\perp}}^{*}$ such that 
$$
-dd^{c}\log\|F\|^{2}
=
\nu\,\widetilde{\omega}_{M^{\perp}}-\delta_{D}
$$
as currents on $({\mathcal M}_{M^{\perp}}^{*})_{\rm reg}$.
Here $\widetilde{\omega}_{M^{\perp}}$ is the current on
$({\mathcal M}_{M^{\perp}}^{*})_{\rm reg}$ defined as
the trivial extension of the K\"ahler form of the Bergman metric,
and $d^{c}=\frac{1}{4\pi i}(\partial-\bar{\partial})$
for a complex manifold. 
\end{itemize}
The notion of meromorphic automorphic form is defined 
in the same manner. 
\end{definition}

Since ${\mathcal B}_{M}$ is a subvariety with 
codimension greater than or equal to $2$ when $r(M)\leq 17$,
the second condition (ii) follows from the first one (i) 
in this case by the Koecher principle.

\subsection
{The equation satisfied by $\tau_{M}$ on the period domain}
\label{subsect:5.2}
\par
Let ${\mathcal A}_{g}$ denote the Siegel modular variety
of degree $g$, which is the coarse moduli space of principally
polarized Abelian varieties of dimension $g$.
The Petersson norm of a Siegel modular form $S$ on the Siegel
upper half space of degree $g$ is denoted by
$\|S\|^{2}$, which is a $C^{\infty}$ function on
${\mathcal A}_{g}$ in the sense of orbifolds.
If $k$ is the weight of $S$, the $(1,1)$-form
$\omega_{{\mathcal A}_{g}}:=-\frac{1}{k}dd^{c}\log\|S\|^{2}$ 
on ${\mathcal A}_{g}$ in the sense of
orbifolds is the K\"ahler form of the Bergman metric.
\par
As an application of Theorem~\ref{Theorem4.3}, we prove 
the automorphic property of $\tau_{M}$ when $r(M)\geq18$. 
For this, we need an extension of \cite[Sect.\,7]{Yoshikawa04}.

\begin{theorem}\label{Theorem5.2}
Let 
$\varPi_{M}\colon\Omega_{M^{\perp}}^{+}\to
{\mathcal M}_{M^{\perp}}$ be the projection and let
$\tau_{\Omega_{M^{\perp}}^{+}}$ be the
$O^{+}(M^{\perp})$-invariant function on 
$\Omega_{M^{\perp}}^{+}\setminus{\mathcal D}_{M^{\perp}}$ 
defined as
$\tau_{\Omega_{M^{\perp}}^{+}}=\varPi_{M}^{*}\tau_{M}$.
Then $\tau_{\Omega_{M^{\perp}}^{+}}$ lies in
$L^{1}_{\rm loc}(\Omega_{M^{\perp}}^{+})$ 
and satisfies the following equation of currents on 
$\Omega_{M^{\perp}}^{+}$:
\begin{equation}\label{eqn:5.1}
dd^{c}\log\tau_{\Omega_{M^{\perp}}^{+}}
=
\frac{r(M)-6}{4}\,\omega_{M}
+J_{M}^{*}\omega_{{\mathcal A}_{g(M)}}
-\frac{1}{4}\delta_{{\mathcal D}_{M^{\perp}}}.
\end{equation}
\end{theorem}

\begin{pf}
Let $O^{+}(M^{\perp})_{[\eta]}\subset O^{+}(M^{\perp})$ 
be the stabilizer of $[\eta]\in\Omega_{M^{\perp}}^{+}$.
As in \cite{Yoshikawa04}, set
$$
H_{d}^{o}:=\{[\eta]\in H_{d};\,
O^{+}(M^{\perp})_{[\eta]}=\{\pm1,\,\pm s_{d}\}\},
\qquad
{\mathcal D}_{M^{\perp}}^{o}
:=
\bigcup_{d\in\Delta_{M^{\perp}}}H_{d}^{o}
$$
and
$Z_{M^{\perp}}:=
\bigcup_{d\in\Delta_{M^{\perp}}}H_{d}\setminus H_{d}^{o}$.
When $r(M)\leq18$, $Z_{M^{\perp}}$ is an analytic subset of 
$\Omega_{M^{\perp}}^{+}$ with codimension greater than or 
equal to $2$ by \cite[Prop.\,1.9 (2)]{Yoshikawa04}.
By \cite[Sect.\,(7.1)]{Yoshikawa04}, 
$\tau_{\Omega_{M^{\perp}}^{+}}$ lies in
$L^{1}_{\rm loc}(\Omega_{M^{\perp}}^{+}\setminus Z_{M^{\perp}})$
and satisfies the following equation of currents on
$\Omega_{M^{\perp}}^{+}\setminus Z_{M^{\perp}}$:
\begin{equation}\label{eqn:5.2}
dd^{c}\log\tau_{\Omega_{M^{\perp}}^{+}}
=
\frac{r(M)-6}{4}\,\omega_{M}
+J_{M}^{*}\omega_{{\mathcal A}_{g(M)}}
-\frac{1}{4}\delta_{{\mathcal D}_{M^{\perp}}}.
\end{equation}
Since ${\rm codim}\,Z_{M^{\perp}}\geq2$ when $r(M)\leq18$, 
we deduce from \eqref{eqn:5.2} and \cite[p.53, Th.\,1]{Siu74} 
that Eq.\eqref{eqn:5.1} holds in this case. 
We consider the case $r(M)\geq19$. 
Since $\Omega_{M^{\perp}}^{+}$ consists of a unique point
when $r(M)=20$, i.e., 
$M^{\perp}\cong\langle2\rangle\oplus\langle2\rangle$,
the assertion is trivial in this case. It suffices to 
prove \eqref{eqn:5.1} when $r(M)=19$, in which case either
$M^{\perp}\cong\langle2\rangle\oplus\langle2\rangle
\oplus\langle-2\rangle$ or
$M^{\perp}\cong{\Bbb U}\oplus\langle2\rangle$
by \cite[p.1434, Table 1]{Nikulin83}.
\par
Assume $M^{\perp}\cong{\Bbb U}\oplus\langle2\rangle$.
By \cite[Th.\,7.1]{Dolgachev96}, there exist isomorphisms
$\Omega_{M^{\perp}}^{+}\cong{\frak H}$ and
$O^{+}(M^{\perp})\cong SL_{2}({\bf Z})$ such that
the $O^{+}(M^{\perp})$-action on $\Omega_{M^{\perp}}^{+}$ 
is identified with the projective action of 
$SL_{2}({\bf Z})$ on ${\frak H}$. Let 
${\mathcal F}:=\{z\in{\frak H};\,|z|\geq1,\,|\Re z|\leq1/2\}$
be the fundamental domain for the $PSL_{2}({\bf Z})$-action
on ${\frak H}$. For $\tau\in{\frak H}$, 
let $SL_{2}({\bf Z})_{\tau}\subset SL_{2}({\bf Z})$ 
be the stabilizer of $\tau$. 
Let $d\in\Delta_{M^{\perp}}$ and
let $z\in{\mathcal F}$ be the point corresponding 
to $[\eta]\in H_{d}$. Since 
$O^{+}(M^{\perp})_{[\eta]}\supset{\bf Z}_{2}\times{\bf Z}_{2}$
and hence $\#\,O^{+}(M^{\perp})_{[\eta]}\geq4$,
we get $\#PSL_{2}({\bf Z})_{z}\geq2$.
By e.g. \cite{Serre70}, we get 
$z\in\{i,e^{\pi/3},e^{2\pi i/3}\}$.
If $z=e^{\pi/3}$ or $e^{2\pi/3}$, then 
$PSL_{2}({\bf Z})_{z}\cong{\bf Z}_{3}$. In this case,
$SL_{2}({\bf Z})_{z}\cong O^{+}(M^{\perp})_{[\eta]}$ does not 
contain a subgroup of order $4$, which contradicts the fact
$O^{+}(M^{\perp})_{[\eta]}\supset\{\pm1,\pm s_{d}\}=
{\bf Z}_{2}\times{\bf Z}_{2}$. 
Hence we get $z=i$. Since $\#SL_{2}({\bf Z})_{i}=4$, we get
$O^{+}(M^{\perp})_{[\eta]}=\{\pm1,\pm s_{d}\}$. 
This implies $H_{d}^{o}=H_{d}$ and 
$Z_{M}=\emptyset$ when $M\cong{\Bbb U}\oplus\langle2\rangle$. 
This proves \eqref{eqn:5.1} when 
$M\cong{\Bbb U}\oplus\langle2\rangle$.
For the case
$M^{\perp}\cong\langle2\rangle\oplus\langle2\rangle
\oplus\langle-2\rangle$, see Sect.\ref{sect:6}.
\end{pf}

\subsection
{The automorphic property of $\tau_{M}$}
\label{subsect:5.3}
\par
\begin{theorem}\label{Theorem5.3}
There exist an integer $\nu\in{\bf Z}_{>0}$
and an (possibly meromorphic) automorphic form $\Psi_{M}$ 
on $\Omega_{M^{\perp}}^{+}$ for $O^{+}(M^{\perp})$ 
of weight $\nu(r(M)-6)$ and a Siegel modular form 
$S_{M}$ on ${\frak S}_{g(M)}$ of weight $4\nu$ such that 
for every $2$-elementary $K3$ surface $(X,\iota)$ of type $M$, 
$$
\tau_{M}(X,\iota)
=
\|\Psi_{M}(\overline{\varpi}_{M}(X,\iota))\|^{-1/2\nu}
\|S_{M}(\varOmega(X^{\iota}))\|^{-1/2\nu}.
$$
Here $\overline{\varpi}_{M}(X,\iota)\in{\mathcal M}_{M^{\perp}}$
denotes the period of $(X,\iota)$ and
$\varOmega(X^{\iota})\in{\mathcal A}_{g(M)}$ denotes
the period of $X^{\iota}$.
\end{theorem}

\begin{pf}
Since the assertion was proved when $r(M)\leq17$ 
(cf.~\cite{Yoshikawa04}), we assume $r(M)\geq18$.
Let $\ell\in{\bf Z}_{>0}$ be sufficiently large.
Let $S$ be a Siegel modular form of weight $4\ell$ on 
${\frak S}_{g(M)}$ such that the function
${\mathcal M}_{M^{\perp}}^{o}\ni(X,\iota)\to
\|S(\varOmega(X^{\iota}))\|^{2}
\in{\bf R}_{\geq0}$ 
does not vanish identically.
Let $F$ be a non-zero automorphic form on 
$\Omega_{M^{\perp}}^{+}$ for $O^{+}(M^{\perp})$ of weight 
$\ell(r(M)-6)$. 
Let $J_{M}^{*}\omega_{{\mathcal A}_{g(M)}}$ be the current
defined as the trivial extension of 
$(J_{M}^{o})^{*}\omega_{{\mathcal A}_{g(M)}}$
from $\Omega_{M^{\perp}}\setminus{\mathcal D}_{M^{\perp}}$
to $\Omega_{M^{\perp}}$, where
$J_{M}^{o}\colon
\Omega_{M^{\perp}}\setminus{\mathcal D}_{M^{\perp}}
\to{\mathcal A}_{g(M)}$ is the holomorphic map defined as
$J_{M}^{o}(\overline{\varpi}_{M}(X,\iota))=\varOmega(X^{\iota})$
(cf. \cite[Sects.\,3.1-3.4]{Yoshikawa04}).
Then the following equations of currents
on $\Omega_{M^{\perp}}^{+}$ hold:
\begin{equation}\label{eqn:5.3}
-dd^{c}\log\|S\|^{2}
=
4\ell\,J_{M}^{*}\omega_{{\mathcal A}_{g(M)}}
-\delta_{J_{M}^{*}{\rm div}(S)},
\end{equation}
\begin{equation}\label{eqn:5.4}
-dd^{c}\log\|F\|^{2}
=
\ell(r(M)-6)\,\omega_{M}-\delta_{{\rm div}(F)}.
\end{equation}
We set
$$
\varphi:=
\tau_{\Omega_{M^{\perp}}^{+}}(\|F\|\cdot\|S\|)^{1/2\ell}.
$$
By \eqref{eqn:5.1}, \eqref{eqn:5.3}, \eqref{eqn:5.4}, 
there is an $O^{+}(M^{\perp})$-invariant ${\bf Q}$-divisor $D$
on $\Omega_{M^{\perp}}^{+}$ satisfying the following equation 
of currents on $\Omega_{M^{\perp}}^{+}$:
\begin{equation}\label{eqn:5.5}
-dd^{c}\log\varphi=\delta_{D}.
\end{equation}
\par
Let $[\eta_{0}]\in\Omega_{M^{\perp}}^{+}$ and 
let $m\in{\bf Z}_{>0}$ be an integer such that $mD$ is
an integral divisor on $\Omega_{M^{\perp}}^{+}$.
We define
$G([\eta]):=
\exp\left(m\int_{[\eta_{0}]}^{[\eta]}\partial\log\varphi\right)$.
Since the residues of the logarithmic $1$-form 
$m\,\partial\varphi$ on $\Omega_{M^{\perp}}^{+}$ are integral, 
$G$ is a meromorphic function on $\Omega_{M^{\perp}}^{+}$ 
with ${\rm div}(G)=m\,D$. By the definition of $G$ and the equality
$\overline{\partial\log\varphi}=\bar{\partial}\log\varphi$, 
we get
\begin{equation}\label{eqn:5.6}
|G([\eta])|^{2}
=
\exp
\left(m\int_{[\eta_{0}]}^{[\eta]}
\partial\log\varphi+\bar{\partial}\log\varphi
\right)
=
\varphi([\eta])^{m}\varphi([\eta_{0}])^{-m}.
\end{equation}
\par
Let $\gamma\in O^{+}(M^{\perp})$.
By the $O^{+}(M^{\perp})$-invariance of $\varphi$, we get
$\gamma^{*}\partial\log\varphi=\partial\log\varphi$, 
which yields that $d\log(\gamma^{*}G/G)=0$. 
Hence there exists a constant
$\chi(\gamma)\in{\bf C}^{*}$ with
\begin{equation}\label{eqn:5.7}
\gamma^{*}G=\chi(\gamma)\,G.
\end{equation}
Since $(\gamma\gamma')^{*}=(\gamma')^{*}\gamma^{*}$ for
$\gamma,\gamma'\in O^{+}(M^{\perp})$, we deduce from 
\eqref{eqn:5.7} that
$\chi\colon O^{+}(M^{\perp})\to{\bf C}^{*}$ is a character.
We see that $|\chi(\gamma)|=1$. Indeed,
by the definition of $\chi$, we get
\begin{equation}\label{eqn:5.8}
\begin{aligned}
|\chi(\gamma)|^{2}
&=
\frac{G(\gamma\cdot[\eta])}{G([\eta])}
\cdot
\overline{\left(\frac{G(\gamma\cdot[\eta])}{G([\eta])}\right)}
\\
&=
\exp\left(m\int_{[\eta_{0}]}^{\gamma\cdot[\eta_{0}]}
\partial\log\varphi\right)\cdot
\exp\left(m\int_{[\eta_{0}]}^{\gamma\cdot[\eta_{0}]}
\bar{\partial}\log\varphi\right)
\\
&=
\exp\left(m\int_{[\eta_{0}]}^{\gamma\cdot[\eta_{0}]}
d\log\varphi\right)
\\
&=
\exp\left[m\,\log\varphi(\gamma\cdot[\eta_{0}])-
m\,\log\varphi([\eta_{0}])\right]
=1,
\end{aligned}
\end{equation}
where the second equality follows from the fact
$\overline{\partial\log\varphi}=\bar{\partial}\log\varphi$
and the last equality follows from the 
$O^{+}(M^{\perp})$-invariance of $\varphi$.
By \eqref{eqn:5.7} and \eqref{eqn:5.8}, 
$G^{-4\ell}F^{m}$ satisfies Definition~\ref{Definition5.1} (1).
\par
Set $C:=\log\varphi([\eta_{0}])$.
By the definition of $\varphi$ and \eqref{eqn:5.6}, we get
\begin{equation}\label{eqn:5.9}
\tau_{\Omega_{M^{\perp}}^{+}}
=
e^{C}|G|^{2/m}(\|F\|\cdot\|S\|)^{-1/2\ell}
=
e^{C}
(\|G^{-4\ell}F^{m}\|^{2}\cdot\|S^{m}\|^{2})^{-1/4m\ell}.
\end{equation}
We set $\nu:=m\ell$, $\Psi_{M}:=G^{-4\ell}F^{m}$ and 
$S_{M}:=S^{m}$. Then
$$
\tau_{M}=(\|\Psi_{M}\|^{2}\|S_{M}\|^{2})^{-1/4\nu}.
$$
Since $S_{M}$ is a Siegel modular form of weight $4\ell m=4\nu$
and since $\Psi_{M}$ is a meromorphic function on
$\Omega_{M^{\perp}}^{+}$ satisfying the functional equation
in Definition~\ref{Definition5.1} (1) with weight 
$(r(M)-6)\ell m=(r(M)-6)\nu$,
it suffices to prove that $\Psi_{M}$ satisfies the regularity
condition (2) in Definition~\ref{Definition5.1}.
Since $|G|$ is an $O^{+}(M^{\perp})$-invariant function
on $\Omega_{M^{\perp}}^{+}\setminus{\mathcal D}_{M^{\perp}}$
by \eqref{eqn:5.6}, we regard $|G|$ as a function on 
${\mathcal M}_{M^{\perp}}^{o}$. Since $F$ is an automorphic 
form on $\Omega_{M^{\perp}}^{+}$ and hence satisfies 
the regularity condition (2) in Definition~\ref{Definition5.1},
it suffices to prove
$\log|G|^{2}\in L^{1}_{\rm loc}(
({\mathcal M}_{M^{\perp}}^{*})_{\rm reg})$ and 
the existence of a ${\bf Q}$-divisor ${\frak D}$ on
${\mathcal M}_{M^{\perp}}^{*}$ satisfying the following equation
of currents on $({\mathcal M}_{M^{\perp}}^{*})_{\rm reg}$:
\begin{equation}\label{eqn:5.10}
-dd^{c}\log|G|^{2}=\delta_{\frak D}.
\end{equation}
\par
Let $\Phi$ be the function on ${\mathcal M}_{M^{\perp}}^{o}$
such that $\varphi=\varPi_{M}^{*}\Phi$.
Let $\overline{D}$ be the closure of $\varPi_{M}(D)$
in ${\mathcal M}_{M^{\perp}}^{*}$. 
By \eqref{eqn:5.5}, we have the equation of currents on 
$({\mathcal M}_{M^{\perp}})_{\rm reg}$:
\begin{equation}\label{eqn:5.11}
-dd^{c}\log\Phi=\delta_{\overline{D}}.
\end{equation}
\par
Let 
${\mathcal B}_{M}=\bigcup_{\alpha\in A}{\mathcal B}_{M,\alpha}$
be the irreducible decomposition. Since $r(M)\geq18$, 
${\mathcal B}_{M}$ is a divisor on ${\mathcal M}_{M^{\perp}}^{*}$.
Let $C\subset{\mathcal M}_{M^{\perp}}^{*}$ be an arbitrary 
irreducible projective curve intersecting ${\mathcal B}_{M}$ 
properly. Let ${\frak b}\in C\cap{\mathcal B}_{M}$ be 
an arbitrary point and let 
$C_{\frak b}=\bigcup_{i\in I}C^{(i)}_{\frak b}$ 
be the irreducible decomposition of the set germ 
$C_{\frak b}=(C,{\frak b})$.
Let $\nu^{(i)}\colon(\varDelta,0)\to C^{(i)}_{\frak b}$ 
be the normalization. By Theorem~\ref{Theorem4.3}, there exists 
$\alpha^{(i)}_{\frak b}\in{\bf Q}$ such that as $t\to 0$,
\begin{equation}\label{eqn:5.12}
\log\tau_{M}|_{C^{(i)}_{\frak b}}(\nu^{(i)}(t))
=
\alpha^{(i)}_{\frak b}\,\log|t|+O(\log(-\log|t|)).
\end{equation}
Since $F$ and $S$ are automorphic forms on $\Omega_{M^{\perp}}^{+}$
and ${\frak S}_{g(M)}$ respectively, there exists 
$\beta^{(i)}_{\frak b},\gamma^{(i)}_{\frak b}\in{\bf Z}$ by 
\cite[Th.\,3.1]{Mumford77} (cf. \cite[Prop.\,3.12]{Yoshikawa04}) 
such that as $t\to0$,
\begin{equation}\label{eqn:5.13}
(\log\|F\|)|_{C^{(i)}_{\frak b}}(\nu^{(i)}(t))
=
\beta^{(i)}_{\frak b}\,\log|t|+O(\log(-\log|t|)),
\end{equation}
\begin{equation}\label{eqn:5.14}
(\log\|S\|)|_{C^{(i)}_{\frak b}}(\nu^{(i)}(t))
=
\gamma^{(i)}_{\frak b}\,\log|t|+O(\log(-\log|t|)).
\end{equation}
By \eqref{eqn:5.12}, \eqref{eqn:5.13}, \eqref{eqn:5.14}, 
there exists $\epsilon^{(i)}_{\frak b}\in{\bf Q}$ such that
\begin{equation}\label{eqn:5.15}
\log\Phi|_{C^{(i)}_{\frak b}}(\nu^{(i)}(t))
=
\epsilon^{(i)}_{\frak b}\,\log|t|+O(\log(-\log|t|))
\qquad
(t\to0).
\end{equation}
By \eqref{eqn:5.11}, \eqref{eqn:5.15}, 
there exists $n_{C,\alpha}\in{\bf Q}$ such that 
the following equation of currents on $C$ holds:
\begin{equation}\label{eqn:5.16}
-dd^{c}\log\Phi|_{C}
=
\delta_{\overline{D}\cap C}
+
\sum_{\alpha\in A}n_{C,\alpha}\,
\delta_{{\mathcal B}_{M,\alpha}\cap C}.
\end{equation}
Since $C\subset{\mathcal M}_{M^{\perp}}^{*}$ is arbitrary, 
this implies that $\partial\log\Phi$ is a logarithmic $1$-form 
on $({\mathcal M}_{M^{\perp}}^{*})_{\rm reg}$ and that
$n_{C,\alpha}$ is the residue of $\partial\log\Phi$
along the irreducible divisor ${\mathcal B}_{M,\alpha}$
for sufficiently general $C$. Since $n_{C,\alpha}$ is 
independent of the choice of sufficiently general $C$, 
we write $n_{\alpha}$ for $n_{C,\alpha}$. 
By \eqref{eqn:5.11} and \eqref{eqn:5.16}, we get 
$\Phi\in L^{1}_{\rm loc}({\mathcal M}_{M^{\perp}}^{*})$ 
and the following equation of currents on 
$({\mathcal M}_{M^{\perp}}^{*})_{\rm reg}$:
\begin{equation}\label{eqn:5.17}
-dd^{c}\log\Phi
=
\delta_{\overline{D}}
+
\sum_{\alpha\in A}n_{\alpha}\,\delta_{{\mathcal B}_{M},\alpha}.
\end{equation}
Set
${\frak D}=m(\overline{D}+
\sum_{\alpha\in A}n_{\alpha}\,{\mathcal B}_{M,\alpha})$. 
Since $\varphi=e^{C}|G|^{2/m}$ by \eqref{eqn:5.6}, 
we get \eqref{eqn:5.10} from \eqref{eqn:5.17}. 
This completes the proof.
\end{pf}

\begin{remark}\label{remark5.4}
In fact, one can prove that the boundary divisor 
${\mathcal B}_{M}$ is irreducible when $r(M)\geq18$. 
The irreducibility of ${\mathcal B}_{M}$ plays a crucial 
role to give an explicit formulae for
$\Psi_{M}$ and $S_{M}$, when $r(M)\geq18$. 
See \cite[Sect.\,11.4]{Yoshikawa09a} for the details.
\end{remark}


\section
{\bf The case $M^{\perp}=\langle2\rangle\oplus\langle2\rangle
\oplus\langle-2\rangle$}
\label{sect:6}
\par
Throughout Sect.\ref{sect:6}, we assume 
$$
M^{\perp}=\langle2\rangle\oplus\langle2\rangle
\oplus\langle-2\rangle
$$ 
and prove that \eqref{eqn:5.1} holds in this case.

\subsection
{Preliminaries}
\label{subsect:6.1}
\par 
Since $M^{\perp}=\langle2\rangle\oplus\langle2\rangle
\oplus\langle-2\rangle$, we get the explicit expression:
$$
\Omega_{M^{\perp}}^{+}
=
\{(x:y:z)\in{\bf P}^{2};\,
x^{2}+y^{2}-z^{2}=0,\,|x|^{2}+|y|^{2}-|z|^{2}>0,\,
|x+iy|>|x-iy|\}.
$$
The unit disc $\varDelta=\{z\in{\bf C};\,|z|<1\}$ 
is isomorphic to $\Omega_{M^{\perp}}^{+}$ by the map
\begin{equation}\label{eqn:6.1}
c\colon
\varDelta\ni z\to
\left(\frac{1+z^{2}}{2}:\frac{1-z^{2}}{2i}:z\right)
\in\Omega_{M^{\perp}}^{+}.
\end{equation}
For $\epsilon\in]0,1[$, we set 
$\varDelta(\epsilon):=\{z\in\varDelta;\,|z|<\epsilon\}$ and
$\Omega_{M^{\perp}}^{+}(\epsilon):=c(\varDelta(\epsilon))$.
We also set $\delta:=(0,0,1)\in\Delta_{M^{\perp}}$. Then
$s_{\delta}(x:y:z):=(x:y:-z)$ is the reflection on
$\Omega_{M^{\perp}}^{+}$ associated to $\delta$, and we have
$$
H_{\delta}\cap\Omega_{M^{\perp}}^{+}=\{c(0)\}=\{(1:-i:0)\}.
$$

\begin{lemma}\label{Lemma6.1}
Let $O^{+}(M^{\perp})_{[\eta]}$ be the stabilizer of 
$[\eta]\in\Omega_{M^{\perp}}^{+}$ in $O^{+}(M^{\perp})$.
Then $\#O^{+}(M^{\perp})_{[\eta]}=8$. 
Moreover, the natural projection 
$\varPi_{M}\colon\Omega_{M^{\perp}}^{+}\to
{\mathcal M}_{M^{\perp}}$ has ramification index $4$
at $c(0)\in\Omega_{M^{\perp}}^{+}$.
\end{lemma}

\begin{pf}
Since 
\begin{equation}\label{eqn:6.2}
O^{+}(M^{\perp})_{c(0)}
=
\langle -1_{M^{\perp}}\rangle
\times
\langle s_{\delta}\rangle
\times
\langle\mu\rangle,
\qquad
s_{\delta}
=
\begin{pmatrix}
1&0&0
\\
0&1&0
\\
0&0&-1
\end{pmatrix},
\quad
\mu
:=
\begin{pmatrix}
0&1&0
\\
-1&0&0
\\
0&0&1
\end{pmatrix},
\end{equation}
we get the first assertion.
Since $-1_{M^{\perp}}$, $s_{\delta}$ and $\mu$ act on 
$\varDelta$ as follows under the identification \eqref{eqn:6.1}:
\begin{equation}\label{eqn:6.3}
-1_{M^{\perp}}(z)=z,
\qquad
s_{\delta}(z)=-z,
\qquad
\mu(z)=i\,z,
\end{equation}
we deduce from \eqref{eqn:6.2}, \eqref{eqn:6.3} that the projection 
$\varPi_{M^{\perp}}\colon\Omega_{M^{\perp}}^{+}\to
{\mathcal M}_{M^{\perp}}$ at $c(0)\in\Omega_{M^{\perp}}^{+}$
is identified with the map ${\bf C}\ni z\to z^{4}\in{\bf C}$
at $z=0$.
\end{pf}

\par
We recall the notion of ordinary singular families of 
$2$-elementary $K3$ surfaces.
Let ${\mathcal Z}$ be a smooth complex threefold.
Let $p\colon{\mathcal Z}\to\varDelta$ be a proper 
surjective holomorphic function without critical 
points on ${\mathcal Z}\setminus p^{-1}(0)$.
Let $\iota\colon{\mathcal Z}\to{\mathcal Z}$ be 
a holomorphic involution preserving the fibers of $p$.
We set $Z_{t}=p^{-1}(t)$ and $\iota_{t}=\iota|_{Z_{t}}$
for $t\in\varDelta$.
Then $p\colon({\mathcal Z},\iota)\to\varDelta$ is 
called an {\it ordinary singular family}
of $2$-elementary $K3$ surfaces of type $M$
if $p$ has a unique, non-degenerate critical point 
on $Z_{0}$ and if $(Z_{t},\iota_{t})$ is a $2$-elementary 
$K3$ surface of type $M$ for all $t\in\varDelta^{*}$. 
See \cite[Sects.\,2.2 and 2.3]{Yoshikawa04}
for more about ordinary singular families of
$2$-elementary $K3$ surfaces.

\begin{proposition}\label{Proposition6.2}
There exist $\epsilon\in]0,1[$ and 
an ordinary singular family of $2$-elementary $K3$ surfaces 
$p\colon({\mathcal Z},\iota)\to
\Omega_{M^{\perp}}^{+}(\epsilon)/\langle s_{\delta}\rangle$ 
of type $M$ with the following properties:
\begin{itemize}
\item[(1)]
The period map for 
$p\colon({\mathcal Z},\iota)\to
\Omega_{M^{\perp}}^{+}(\epsilon)/\langle s_{\delta}\rangle$
is given by the projection 
$\Omega_{M^{\perp}}^{+}(\epsilon)/\langle s_{\delta}\rangle
\to{\mathcal M}_{M^{\perp}}$.
\item[(2)]
The map $p\colon{\mathcal Z}\to
\Omega_{M^{\perp}}^{+}(\epsilon)/\langle s_{\delta}\rangle$ 
is projective.
\end{itemize}
\end{proposition}

\begin{pf}
We follow \cite[Th.\,2.6]{Yoshikawa04}. 
By Lemma~\ref{Lemma6.1}, we get $H_{\delta}^{o}=\emptyset$.
In particular, \cite[Th.\,2.6]{Yoshikawa04} does not
apply at once for
$M^{\perp}=\langle2\rangle\oplus\langle2\rangle
\oplus\langle-2\rangle$. 
However, in the proof of \cite[Th.\,2.6]{Yoshikawa04}, 
the fact $c(0)\in H_{\delta}^{o}$ was used only to 
deduce the following (i), (ii):
\begin{itemize}
\item[(i)]
$s_{\delta}(c(t))=c(-t)$ for all $t\in\varDelta(\epsilon)$.
\item[(ii)]
Under the inclusion 
$M^{\perp}=\langle2\rangle\oplus\langle2\rangle
\oplus\langle-2\rangle\subset{\Bbb L}_{K3}$, set 
$$
\Delta_{c(0)}:=\{d\in\Delta_{{\Bbb L}_{K3}};\,
\langle d,(1,-i,0)\rangle=0\}.
$$ 
Then there exists $m\in M$ such that 
$\langle m,m\rangle_{M}>0$ and 
$m^{\perp}\cap\Delta_{c(0)}=\{\pm\delta\}$.
\end{itemize}
Once (i), (ii) are verified, the proof of
\cite[Th.\,2.6]{Yoshikawa04} for the existence of 
an ordinary singular family of $2$-elementary $K3$ surfaces 
$p\colon({\mathcal Z},\iota)\to
\Omega_{M^{\perp}}^{+}(\epsilon)/\langle s_{\delta}\rangle$ 
with (1), (2) works. (Notice that the condition $r(M)\leq17$ 
was not used in \cite[Th.\,2.6]{Yoshikawa04}.) 
Hence it suffices to
prove (i), (ii).
By \eqref{eqn:6.1}, we get (i).
By \cite[Lemma A.2]{Yoshikawa04}, it suffices to prove
$M^{\perp}\cap\Delta_{c(0)}=\{\pm\delta\}$.
Since $c(0)=(1:-i:0)$, we get
$M^{\perp}\cap\Delta_{c(0)}=\{\pm\delta\}$.
This proves (ii). 
\end{pf}

\begin{proposition}\label{Proposition6.3}
Set 
${\frak q}:=\varPi_{M^{\perp}}(c(0))\in{\mathcal M}_{M^{\perp}}$.
Then there exist a pointed smooth projective curve 
$(B,{\frak p})$, a neighborhood $U$ of ${\frak p}$, 
a holomorphic map between curves
$\varphi\colon(B,{\frak p})\to
({\mathcal M}_{M^{\perp}}^{*},{\frak q})$,
a smooth projective threefold ${\mathcal X}$ with an involution 
$\theta\colon{\mathcal X}\to{\mathcal X}$, 
and a surjective holomorphic map $p\colon{\mathcal X}\to B$
with the following properties:
\begin{itemize}
\item[(1)]
$\varphi(B)={\mathcal M}_{M^{\perp}}^{*}$ and
the map $f|_{U}\colon(U,{\frak p})\to(\varphi(U),{\frak q})$ 
is a double covering with a unique ramification 
point ${\frak p}$.
\item[(2)]
The map $p\colon{\mathcal X}\to B$ is ${\bf Z}_{2}$-equivariant 
with respect to the ${\bf Z}_{2}$-action on ${\mathcal X}$ 
induced by $\theta$ and with respect to the trivial 
${\bf Z}_{2}$-action on $B$.
\item[(3)]
The family of algebraic surfaces with involution
$p|_{p^{-1}(U)}\colon({\mathcal X},\theta)|_{p^{-1}(U)}\to U$ 
is an ordinary singular family of $2$-elementary $K3$ surfaces 
of type $M$ with period map $\varphi|_{U}$.
\end{itemize}
\end{proposition}

\begin{pf}
We follow \cite[Th.\,2.8]{Yoshikawa04}. 
Since $H_{\delta}^{o}=\emptyset$ and hence
${\mathcal D}_{M^{\perp}}^{o}=\emptyset$ by 
Lemma~\ref{Lemma6.4} below,
\cite[Th.\,2.8]{Yoshikawa04} does not apply at once for
$M^{\perp}=\langle2\rangle\oplus\langle2\rangle
\oplus\langle-2\rangle$. 
Set 
$S:=\Omega_{M^{\perp}}^{+}(\epsilon)/\langle s_{\delta}\rangle
\cong\varDelta$. Let $\bar{c}(0)\in S$ be the image of $c(0)$
and let
$\gamma\colon(S,\bar{c}(0))\to({\mathcal M}_{M^{\perp}},{\frak p})$ 
be the projection induced from $\varPi_{M^{\perp}}$. 
By Proposition~\ref{Proposition6.2}, there is an ordinary
singular family of $2$-elementary $K3$ surfaces
$p\colon({\mathcal Z},\iota)\to S$ of type $M$ with 
period map $\gamma$.
We set $C:={\mathcal M}_{M^{\perp}}^{*}$.
By Theorem~\ref{Theorem3.1} (2) applied to 
$p\colon({\mathcal Z},\iota)\to S$,
there exist $\varphi\colon B\to C$, $f\colon X\to B$ and
$\theta\colon X\to X$ as in Theorem~\ref{Theorem3.1} satisfying 
(a), (b), (c), (d), (e), (f). 
\par
We prove that
${\rm Sing}\, X\cap X_{\frak p}=\emptyset$.
Since $(X_{\frak p},\theta|_{X_{\frak p}})\cong
(Z_{\bar{c}(0)},\iota|_{Z_{\bar{c}(0)}})$ by 
Theorem~\ref{Theorem3.1} (d) and since $Z_{0}$ has a unique 
$A_{1}$-singularity $o:={\rm Sing}\,Z_{0}$, the deformations
$p\colon(X,X_{\frak p})\to(B,{\frak p})$ and 
$p\colon({\mathcal Z},Z_{0})\to(S,\bar{c}(0))$ induce maps 
$\rho_{f}\colon(B,{\frak p})\to{\rm Def}(A_{1})$ 
and $\rho_{p}\colon(\varDelta,0)\to{\rm Def}(A_{1})$,
where ${\rm Def}(A_{1})\cong({\bf C},0)$ is the Kuranishi space
of $2$-dimensional $A_{1}$-singularity. Since $(Z_{0},o)$
is an $A_{1}$-singularity, there is an isomorphism of germs 
$(X,o)\cong({\mathcal Z},o)$ if and only if
${\rm mult}_{\frak p}\rho_{f}={\rm mult}_{0}\rho_{p}$.
\par
Recall that
${\rm Def}(Z_{0})$ is the Kuranishi space of $Z_{0}$ and
let $\rho\colon{\rm Def}(Z_{0})\to{\rm Def}(A_{1})$ be the map
of germs induced by the local semiuniversal deformation of
$Z_{0}$ over ${\rm Def}(Z_{0})$.
Recall that $\mu_{f}\colon(B,{\frak p})\to{\rm Def}(Z_{0})$ 
and $\mu_{p}\colon(\varDelta,0)\to{\rm Def}(Z_{0})$
are the maps induced by the deformations
$f\colon(X,X_{\frak p})\to(B,{\frak p})$ and 
$p\colon({\mathcal Z},Z_{0})\to(\varDelta,0)$, respectively.
Since 
$\rho_{f}=\rho\circ\mu_{f}$ and $\rho_{p}=\rho\circ\mu_{p}$,
there exists by Theorem~\ref{Theorem3.1} (f) an isomorphism 
of germs $\psi\colon(\varDelta,0)\cong(B,{\frak p})$ such that
$\rho_{f}\circ\psi(t)-\rho_{p}(t)\in t^{\nu}{\bf C}\{t\}$.
By choosing $\nu\geq2$ in Theorem~\ref{Theorem3.1} (2), 
this implies that
${\rm mult}_{\frak p}\rho_{f}={\rm mult}_{0}\rho_{p}$.
Since ${\mathcal Z}$ is smooth, we get
${\rm Sing}\, X\cap X_{\frak p}=\emptyset$.
Let $U$ be a small neighborhood of ${\frak p}$ in $B$.
\par
Since the map 
$\gamma\colon(S,c(0))\to(C,{\frak q})$ has ramification 
index $2$ by Lemma~\ref{Lemma6.1}, we get (1) by 
Theorem~\ref{Theorem3.1} (e).
We get (2) by Theorem~\ref{Theorem3.1} (b). 
By Theorem~\ref{Theorem3.1} (c),
$(X_{b},\theta|_{X_{b}})$ is a $2$-elementary $K3$ surface 
of type $M$ with period $\varphi(b)$ for 
$b\in U\setminus\{{\frak p}\}$. 
Since $p^{-1}(U)$ is smooth,
$p|_{p^{-1}(U)}\colon(X,\theta)|_{p^{-1}(U)}\to U$
is an ordinary singular family of $2$-elementary $K3$ surfaces 
of type $M$ with period map $\varphi|_{U}$. This proves (3).
\par
There is a resolution ${\mathcal X}\to X$ such that 
$\theta$ lifts to an involution
$\widetilde{\theta}\colon{\mathcal X}\to{\mathcal X}$
\cite[Th.\,13.4]{BierstoneMilman97}.
Replace $(X,\theta)$ by $({\mathcal X},\widetilde{\theta})$.
Then (1), (2), (3) are satisfied and ${\mathcal X}$ is smooth.
This completes the proof.
\end{pf}

\begin{lemma}\label{Lemma6.4}
The group $O^{+}(M^{\perp})$ acts transitively
on $\Delta_{M^{\perp}}/\pm1$. 
\end{lemma}

\begin{pf}
Let $L$ be an odd unimodular lattice of signature $(2,1)$
and let $\delta,\delta'\in L$ be vectors with
$\langle\delta,\delta\rangle=\langle\delta',\delta'\rangle=-1$.
Set $\Lambda:=\delta^{\perp}$ and $\Lambda':=(\delta')^{\perp}$.
Since $\Lambda^{\lor}/\Lambda\cong
({\bf Z}\delta)^{\lor}/({\bf Z}\delta)$ and
$(\Lambda')^{\lor}/\Lambda'\cong
({\bf Z}\delta')^{\lor}/({\bf Z}\delta')$, the equality 
$\langle\delta,\delta\rangle=\langle\delta',\delta'\rangle=-1$
implies that $\Lambda$ and $\Lambda'$ are positive-definite
unimodular lattices of rank $2$. It is classical that
$\Lambda\cong\Lambda'$. 
Since ${\bf Z}\delta\oplus\Lambda\subset L$ and since
both of ${\bf Z}\delta\oplus\Lambda$ and $L$ are unimodular,
we get $L={\bf Z}\delta\oplus\Lambda$. Similarly, 
$L={\bf Z}\delta'\oplus\Lambda'$. 
Let $\varphi\colon\Lambda\to\Lambda'$ be an isometry. Then
$g_{\pm}\colon{\bf Z}\delta\oplus\Lambda\ni m\delta+\lambda\to
\pm m\delta'+\varphi(\lambda)\in{\bf Z}\delta'\oplus\Lambda'$
is an isometry of $L$ with $g_{\pm}(\delta)=\pm\delta'$ such that
either $g_{+}\in O^{+}(L)$ or $g_{-}\in O^{+}(L)$.
This proves the lemma.
\end{pf}

\subsection
{Proof of (5.1) }
\label{subsect:6.2}
\par
By \cite[Th.\,5.9]{Yoshikawa04}, we have
the following equation of $C^{\infty}$ $(1,1)$-forms on 
$\Omega_{M^{\perp}}^{+}\setminus{\mathcal D}_{M^{\perp}}$
\begin{equation}\label{eqn:6.4}
dd^{c}\log\tau_{\Omega_{M^{\perp}}^{+}}
=
\frac{r(M)-6}{4}\,\omega_{M}
+J_{M}^{*}\omega_{{\mathcal A}_{g(M)}}.
\end{equation}
By Theorem~\ref{Theorem4.3} and \eqref{eqn:6.4}, there exists
$m_{d}\in{\bf Q}$ for every $d\in\Delta_{M^{\perp}}$ 
such that the following equation of currents on
$\Omega_{M^{\perp}}^{+}$ holds:
\begin{equation}\label{eqn:6.5}
dd^{c}\log\tau_{\Omega_{M^{\perp}}^{+}}
=
\frac{r(M)-6}{4}\,\omega_{M}
+
J_{M}^{*}\omega_{{\mathcal A}_{g(M)}}
-
\sum_{d\in\Delta_{M^{\perp}}/\pm}m(d)\,\delta_{H_{d}}.
\end{equation}
\par
We compute $m(\delta)$ for 
$\delta=(0,0,1)\in\Delta_{M^{\perp}}$. 
In Proposition~\ref{Proposition6.3}, we may assume that $U$ is 
equipped with a coordinate function $u$ centered at ${\frak p}$.
By \cite[Th.\,7.5]{Yoshikawa04} applied to the ordinary
singular family 
$p|_{p^{-1}(U)}\colon(X,\theta)|_{p^{-1}(U)}\to U$
in Proposition~\ref{Proposition6.3} (3), we get
\begin{equation}\label{eqn:6.6}
\log\tau_{M}(X_{u},\theta|_{X_{u}})
=
-\frac{1}{8}\log|u|^{2}+O\left(\log(-\log|u|^{2})\right)
\qquad
(u\to0).
\end{equation}
Let $t$ be a coordinate function on $f(U)$ centered at ${\frak p}$. 
By Proposition~\ref{Proposition6.3} (1), there exists 
$\epsilon(u)\in{\mathcal O}(U)$ with $\epsilon(0)\not=0$
such that
\begin{equation}\label{eqn:6.7}
t\circ f(u)=u^{2}\epsilon(u).
\end{equation}
By \eqref{eqn:6.6} and \eqref{eqn:6.7}, we get
\begin{equation}\label{eqn:6.8}
\begin{aligned}
\log\tau_{M}(f(u))
&=
\log\tau_{M}(X_{u},\theta|_{X_{u}})
\\
&=
-\frac{1}{8}\log|u|^{2}+O\left(\log(-\log|u|^{2})\right)
\\
&=
-\frac{1}{16}\log|t\circ f(u)|^{2}
+O\left(\log(-\log|t\circ f(u)|^{2})\right).
\end{aligned}
\end{equation}
Since the projection
$\varPi_{M^{\perp}}\colon\Omega_{M^{\perp}}^{+}
\to{\mathcal M}_{M}$ has ramification index $4$ at $c(0)$
by Lemma~\ref{Lemma6.1}, we get by \eqref{eqn:6.8}
\begin{equation}\label{eqn:6.9}
\begin{aligned}
\tau_{\Omega_{M^{\perp}}^{+}}(c(z))
&=
-\frac{1}{16}({\rm index}_{c(0)}\varPi_{M^{\perp}})\,\log|z|^{2}
+O\left(\log(-\log|z|^{2})\right)
\\
&=
-\frac{1}{4}\log|z|^{2}
+O\left(\log(-\log|z|^{2})\right)
\qquad
(z\to0).
\end{aligned}
\end{equation}
By \eqref{eqn:6.9}, we get $m(\delta)=\frac{1}{4}$.
Since $\Delta_{M^{\perp}}/\pm1=O^{+}(M^{\perp})\cdot\delta$
by Lemma~\ref{Lemma6.4} and since $\tau_{\Omega_{M^{\perp}}^{+}}$ 
is $O^{+}(M^{\perp})$-invariant, we get
$m(d)=m(\delta)=\frac{1}{4}$ for all $d\in\Delta_{M^{\perp}}$.
Substituting $m(d)=\frac{1}{4}$ into \eqref{eqn:6.5}, 
we get \eqref{eqn:5.1}.
This completes the proof of \eqref{eqn:5.1}.
\qed



\begin{thebibliography}{99}

\bibitem{BierstoneMilman97}
Bierstone, E., Milman, P.
\newblock
{\em Canonical desingularization in characteristic zero
by blowing up the maximum strata of a local invariant},
\newblock 
Invent. Math.
\newblock
{\bf 128}
\newblock 
(1997),
\newblock 
207--302.


\bibitem{Bismut95}
Bismut, J.-M.
\newblock 
{\em Equivariant immersions and Quillen metrics},
\newblock 
J. Differential Geom.
\newblock 
{\bf 41}
\newblock 
(1995),
\newblock 
53--157.


\bibitem{BGS88}
Bismut, J.-M., Gillet, H., Soul\'e, C.
\newblock 
{\em Analytic torsion and 
holomorphic determinant bundles I,II,III},
\newblock 
Commun. Math. Phys.
\newblock 
{\bf 115}
\newblock 
(1988),
\newblock 
49-78, 
\newblock 
79-126, 
\newblock 
301-351.


\bibitem{DLS94}
Demailly, J.-P., Lempert, L., Shiffman, B.
\newblock 
{\em Algebraic approximations of holomorphic maps from
Stein domains to projective manifolds},
\newblock 
Duke Math. J.
\newblock
{\bf 76}
\newblock 
(1994),
\newblock 
333--363.


\bibitem{Dolgachev96}
Dolgachev, I.
\newblock
{\em Mirror symmetry for lattice polarized $K3$ surfaces},
\newblock
J. Math. Sci.
\newblock
{\bf 81}
\newblock
(1996),
\newblock
2599--2630.


\bibitem{FinashinKharlamov08}
Finashin, S., Kharlamov, V.
\newblock
{\em Deformation classes of real four-dimensional
cubic hypersurfaces},
\newblock
J. Algebraic Geom.,
\newblock
{\bf 17}
\newblock
(2008),
\newblock
677--707.


\bibitem{KohlerRoessler01}
K\"ohler, K., Roessler, D.
\newblock
{\em A fixed point formula of Lefschetz type in 
Arakelov geometry I},
\newblock Invent. Math.
\newblock
{\bf 145}
\newblock 
(2001),
\newblock 
333--396.


\bibitem{Ma00}
Ma, X.
\newblock 
{\em Submersions and equivariant Quillen metrics},
\newblock 
Ann. Inst. Fourier
\newblock 
{\bf 50}
\newblock 
(2000),
\newblock 
1539-1588.


\bibitem{Mumford77}
Mumford, D.
\newblock
{\em Hirzebruch's proportionality theorem in the
non-compact case},
\newblock Invent. Math.
\newblock 
{\bf 42}
\newblock 
(1977),
\newblock 
239--272.


\bibitem{Nikulin80}
Nikulin, V.V.
\newblock 
{\em Integral Symmetric bilinear forms
and some of their applications},
\newblock 
Math. USSR Izv.
\newblock 
{\bf 14}
\newblock 
(1980),
\newblock 
103--167.


\bibitem{Nikulin83}
\bysame
\newblock 
{\em Factor groups of groups of automorphisms of 
hyperbolic forms with respect to subgroups 
generated by $2$-reflections},
\newblock 
J. Soviet Math.
\newblock 
{\bf 22}
\newblock 
(1983),
\newblock 
1401--1476.


\bibitem{RaySinger73}
Ray, D.B., Singer, I.M.
\newblock
{\em Analytic torsion for complex manifolds},
\newblock 
Ann. Math.
\newblock
{\bf 98}
\newblock 
(1973),
\newblock 
154--177.


\bibitem{Serre70}
Serre, J.-P.
\newblock 
{\em Cours d'arithm\'etique},
\newblock 
Presses Universitaires de France,
\newblock 
Paris
\newblock 
(1970)


\bibitem{Siu74}
Siu, Y.-T.
\newblock {\em Analyticity of sets associated to Lelong 
numbers and the extension of closed positive currents},
\newblock 
Invent. Math.
\newblock 
{\bf 27}
\newblock 
(1974),
\newblock 
53--156.


\bibitem{Takegoshi95}
Takegoshi, K.
\newblock {\em Higher direct images of canonical
sheaves tensorized with semi-positive vector bundles
by proper K\"ahler morphisms},
\newblock Math. Ann.
\newblock {\bf 303}
\newblock (1995),
\newblock 389--416


\bibitem{Yoshikawa04}
Yoshikawa, K.-I.
\newblock
{\em $K3$ surfaces with involution,
equivariant analytic torsion,
and automorphic forms on the moduli space},
\newblock
Invent. Math.
\newblock
{\bf 156}
\newblock
(2004),
\newblock
53-117.


\bibitem{Yoshikawa07}
\bysame
\newblock 
{\em On the singularity of Quillen metrics},
\newblock 
Math. Ann. 
\newblock 
{\bf 337}
\newblock 
(2007),
\newblock 
61--89.


\bibitem{Yoshikawa09a}
\bysame
\newblock
{\em $K3$ surfaces with involution,
equivariant analytic torsion,
and automorphic forms on the moduli space II},
\newblock
preprint
\newblock
(2009)


\bibitem{Yoshikawa09b}
\bysame
\newblock 
{\em Singularities and analytic torsion},
\newblock 
preprint
\newblock 
(2010)

\end{thebibliography}
\end{document}